\documentclass[a4paper,12pt,reqno]{amsart}
\usepackage[a4paper,margin=1.9cm,top=2.5cm,bottom=2.5cm,centering,vcentering]{geometry}
\usepackage{amssymb,mathtools,cite,enumerate,color,eqnarray,hyperref}
\definecolor{blue}{RGB}{0, 0, 200}
\definecolor{pink}{RGB}{252, 0, 50}

\hypersetup{
           colorlinks=true,
           linkcolor=blue,
           urlcolor=red,
           citecolor=pink,
          }

\theoremstyle{plain}
\newtheorem{theorem}{Theorem}[section]

\newtheorem{lemma}[theorem]{Lemma}

\theoremstyle{definition}

\newtheorem{Conjecture}[theorem]{Conjecture}

\numberwithin{equation}{section}

\numberwithin{equation}{section}

\begin{document}
																			
\title{New congruences for partitions where the even parts are distinct}

\author[H. Nath]{Hemjyoti Nath}
\address{Department of mathematical sciences, Tezpur University, Napaam, Tezpur, Assam 784028, India}
\email{hemjyotinath40@gmail.com}


%\date{\today.}

\begin{abstract}
We denote the number of partitions of $n$ wherein the even parts are distinct (and the odd parts are unrestricted) by $ped(n)$. In this paper, we will use generating function manipulations to obtain new congruences for $ped(n)$ modulo $24$.
\end{abstract}

\maketitle

\section{Introduction and Main Result}
A partition of a positive integer $n$ is a non-increasing sequence of positive integers whose sum is equal to $n$. If $p(n)$ denotes the number of partitions of a positive integer $n$ and we adopt the convention $p(0)=1$, then the generating function for $p(n)$ satisfies the identity
\begin{equation*}
    \sum_{n=0}^{\infty}p(n)q^n = \frac{1}{(q;q)_{\infty}},
\end{equation*}
where
\begin{equation*}
    (a;q)_{\infty} := \prod_{n=0}^{\infty}(1-aq^n), \quad |q|<1.
\end{equation*}
Throughout this paper, we write
\begin{equation*}
    f_k := (q^k;q^k)_{\infty}, \quad \text{for any integer} \quad k\geq1.
\end{equation*}
The number of partitions of $n$ wherein the even parts are distinct (and the odd parts are unrestricted) is denoted by $ped(n)$. The generating function for $ped(n)$ \cite{1} is 
\begin{equation}\label{e0}
    \sum_{n=0}^{\infty} ped(n)q^n = \frac{(-q^2;q^2)_{\infty}}{(q;q^2)_{\infty}} = \frac{(q^4;q^4)_{\infty}}{(q;q)_{\infty}}.
\end{equation}
Note that by $\eqref{e0}$, the number of partitions of $n$ wherein the even parts are distinct (and the odd parts are unrestricted) equals the number of partitions of $n$ with no parts divisible by $4$, i.e., the $4$-regular partitions (see \cite{1} and references therein). In recent years many congruences for the number of $4$-regular partitions have been discovered (see \cite{8,9,10,11,12,13,14,15} and references therein).

Numerous congruence properties are known for the function $ped(n)$. For example, Andrews, Hirschhorn and Sellers \cite{1} proved that for $\alpha \geq 1$ and $n\geq 0$,
\begin{align*}
    ped(3n+2) & \equiv 0 \pmod{2},\\
    ped(9n+4) & \equiv 0 \pmod{4},\\
    ped(9n+7) & \equiv 0 \pmod{12},\\
    ped\left( 3^{2\alpha+2}n+\frac{11\cdot 3^{2\alpha+1}-1}{8} \right) & \equiv 0 \pmod{2},\\
    ped\left( 3^{2\alpha+1}n+\frac{17\cdot 3^{2\alpha}-1}{8} \right) & \equiv 0 \pmod{6},\\
    ped\left( 3^{2\alpha+2}n+\frac{19\cdot 3^{2\alpha+1}-1}{8} \right) & \equiv 0  \pmod{6}.
\end{align*}
Recently, Xia \cite{7} obtained many interesting infinite families of congruences modulo $8$ for $ped(n)$.

The aim of this paper is to establish new congruences modulo $24$ for $ped(n)$. In the next theorem, we state our main results.

\begin{theorem}\label{t1}
    For every $n\geq0$, we have
\begin{align*}
     ped(225n+43) & \equiv 0 \pmod{24},\\
     ped(225n+88) & \equiv 0 \pmod{24},\\
     ped(225n+133) & \equiv 0 \pmod{24},\\
     ped(225n+223) & \equiv 0 \pmod{24}.
\end{align*}
Furthermore, for every $k \geq 1$ and $n \geq 0$, we have
    \begin{equation*}
    ped(9n+7) \equiv ped\left( 9\cdot5^{2k}n + \frac{57\cdot5^{2k}-1}{8}\right) \pmod{24}.
\end{equation*}
\end{theorem}

The paper is organised as follows: In Section \ref{sec:pre}, we present some preliminaries required for our proofs. In Sections \ref{sec:t1}, we  present the proof of Theorem \ref{t1}.

\section{Preliminaries}\label{sec:pre}
In this section, we collect the $q$-series identities that are used in our proofs. Recall that Ramanujan's general theta function $f(a,b)$ \cite{6} is defined by
\begin{equation*}
    f(a,b)=\sum_{n=-\infty}^{\infty}a^{n(n+1)/2}b^{n(n-1)/2}, \quad |ab|<1.
\end{equation*}
Important special cases of $f(a,b)$ are the theta functions $\varphi(q)$, $\psi(q)$ and $f(-q)$, which satisfies the identities
\begin{equation*}
    \varphi(q) := f(q,q) = \sum_{n=-\infty}^{\infty}q^{n^2} = (-q;q^2)_{\infty}^2(q^2;q^2)_{\infty} = \frac{f_2^5}{f_1^2f_4^2},
\end{equation*}
\begin{equation*}
    \psi(q) := f(q,q^3) = \sum_{n=0}^{\infty}q^{n(n+1)/2} = \frac{(q^2;q^2)_{\infty}}{(q;q^2)_{\infty}} = \frac{f_2^2}{f_1},
\end{equation*}
and
\begin{equation*}
    f(-q):=f(-q,-q^2)=\sum_{n=-\infty}^{\infty}(-1)^{n}q^{n(3n+1)/2} = (q;q)_{\infty}=f_1.
\end{equation*}
In terms of $f(a,b)$, Jacobi's triple product identity \cite{6} is given by
\begin{equation*}
    f(a,b) = (-a;ab)_{\infty}(-b;ab)_{\infty}(ab;ab)_{\infty}.
\end{equation*}

\begin{lemma}[Hirschhorn \cite{5}]
    We have that
    \begin{equation}\label{e1}
        f_1 = f_{25}\left(R(q^5)-q-q^2R(q^5)^{-1}\right),
    \end{equation}
where 
\begin{equation*}
    R(q) = \frac{(q;q^5)_{\infty}(q^4;q^5)_{\infty}}{(q^2;q^5)_{\infty}(q^3;q^5)_{\infty}}.
\end{equation*}
\end{lemma}

\section{Proof of theorem $\ref{t1}$}\label{sec:t1}
Andrews, Hirschhorn and Sellers \cite{1} proved that
\begin{equation}\label{e2}
     \sum_{n=0}^{\infty}ped(9n+7)q^n = 12 \frac{f_2^4 f_3^6 f_4}{f_1^{11}}.
\end{equation}
Therefore, 
\begin{equation}\label{e3}
    ped(9n+7) \equiv 0 \pmod{12}.
\end{equation}
It follows from $\eqref{e2}$ that
\begin{align}\label{e3a}
    \sum_{n=0}^{\infty} ped(9n+7)q^n \equiv 12 \frac{f_2^4 f_3^6 f_4}{f_1^{11}} \pmod {24}.
\end{align}
But, by the binomial theorem, $f_t^{2m} \equiv f_{2t}^m \pmod{2}$, for all positive integers $t$ and $m$. 

Therefore, it follows from $\eqref{e3a}$ that 
\begin{equation}\label{e3b}
    \sum_{n=0}^{\infty} ped(9n+7)q^n \equiv 12 f_1f_6f_{12} \pmod {24}.
\end{equation}

Employing $\eqref{e1}$ in $\eqref{e3b}$, we arrive at
\begin{align}\label{e4}
\sum_{n=0}^{\infty} ped(9n+7) q^n & \equiv 12f_{25}f_{150}f_{300} \left(
R_{30} R_5 R_{60} - R_{30} R_{60} q^2 - \frac{R_{30} R_{60} q^2}{R_5} + R_5 R_{60} q^6 \right. \nonumber \\ 
&\quad + R_{60} q^7 + \frac{R_{60} q^8}{R_5}- R_{30} R_5 q^{12} - \frac{R_5 R_{60} q^{12}}{R_{30}} + R_{30} q^{13} + \frac{R_{60} q^{13}}{R_{30}} \nonumber \\ 
&\quad + \frac{R_{30} q^{14}}{R_5} + \frac{R_{60} q^{14}}{R_{30} R_5} + R_5 q^{18} - q^{19} - \frac{q^{20}}{R_5} + \frac{R_5 q^{24}}{R_{30}} - \frac{R_{30} R_5 q^{24}}{R_{60}} \nonumber \\ 
&\quad - \frac{q^{25}}{R_{30}} + \frac{R_{30} q^{25}}{R_{60}} - \frac{q^{26}}{R_{30} R_5} + \frac{R_{30} q^{26}}{R_5 R_{60}} + \frac{R_5 q^{30}}{R_{60}} - \frac{q^{31}}{R_{60}} - \frac{q^{32}}{R_5 R_{60}} \nonumber  \\
&\quad \left.  + \frac{R_5 q^{36}}{R_{30} R_{60}} - \frac{q^{37}}{R_{30} R_{60}} - \frac{q^{38}}{R_{30} R_5 R_{60}}\right) \pmod{24}.
\end{align}
Extracting the terms involving $q^{5n+4}$ from both sides of $\eqref{e4}$, dividing both sides by $q^4$ and then replacing $q^5$ by $q$, yields
\begin{equation*}
     \sum_{n=0}^{\infty} ped(9(5n+4)+7)q^n \equiv 12f_{5}f_{30}f_{60}\left( 2q^2\frac{R_6}{R_1}-q^3 \right)\pmod{24},
\end{equation*}
from which it follows that
\begin{equation}\label{e5}
     \sum_{n=0}^{\infty} ped(45n+43)q^n \equiv 12q^3f_{5}f_{30}f_{60}\pmod{24}.
\end{equation}

Next, equating the coefficients of $q^{5n+j}$ on both sides of this congruence, where $j=0,1,2,4$, gives the congruences in Theorem $\ref{t1}$.

Further, extracting the terms involving $q^{5n+3}$ from both sides of $\eqref{e5}$, dividing both sides by $q^3$ and then replacing $q^5$ by $q$, yields
\begin{equation*}
     ped(225n+178) \equiv 12f_{1}f_{6}f_{12}\pmod{24}.
\end{equation*}
which is equivalent to 
\begin{equation}\label{e6}
     ped(9n+7) \equiv ped(225n+178) \pmod{24}.
\end{equation}
Successive iterations of $\eqref{e6}$ give
\begin{align*}
   ped(9n+7) & \equiv ped(9(25n+19)+7)  \\ 
   & \equiv ped(225(25n+19)+178) \\ 
   & \equiv ped(9\cdot5^4n +9\cdot5^2\cdot19+9\cdot19+7)  \\ 
   & \quad\vdots \\
   & \equiv ped(9\cdot5^{2k} n + 9\cdot19\cdot5^{2k-2} + \ldots + 9\cdot19 +7 ) \\
   & \equiv ped\left( 9\cdot5^{2k}n + \frac{57\cdot5^{2k}-1}{8}\right) \pmod{24}.
\end{align*}
This completes the proof.\\

The author would like to end this section with the following conjecture: 
\begin{Conjecture}
 For each nonnegative integer $n$,
\begin{align*}
     ped(225n+43) & \equiv 0 \pmod{192},\\
     ped(225n+88) & \equiv 0 \pmod{192},\\
     ped(225n+133) & \equiv 0 \pmod{192},\\
     ped(225n+223) & \equiv 0 \pmod{192}.
\end{align*}    
\end{Conjecture}

\section{Concluding remarks}
Recently, Chen \cite{8} proved some vanishing results on the coefficients of $\theta_{\chi}(z)$ and the product of two theta functions. Using these results and some generating function manipulations we can find many more congruences for $ped(n)$ modulo $24$.

\end{document}